\documentstyle[12pt,xypic]{amsart}

\numberwithin{equation}{section}
\newtheorem{thm}{Theorem}[section]
\newtheorem{prop}[thm]{Proposition}
\newtheorem{cor}[thm]{Corollary}

\newtheorem{lemma}[thm]{Lemma}

\newenvironment{proof}[1]{\vspace*{2mm}\noindent {\it Proof }{\it
#1.}}{\hspace{2mm}$\Box$\vspace*{2mm}\newline}

\newcommand{\T}{Teichm\"{u}ller }

\newcommand{\C}{{\Bbb C}}
\newcommand{\R}{{\Bbb R}}
\newcommand{\Z}{{\Bbb Z}}

\newcommand{\HS}{{\cal H}_{\Bbb R}^1(\Sigma)}
\newcommand{\HSc}{{\cal H}_{\Bbb R}^1(\Sigma^c)}

\begin{document}
\title[Jacobian varieties of Klein surfaces]
{Torelli groups and Jacobian varieties of non-orientable compact Klein
surfaces}
\author[Pablo Ar\'es Gastesi]{Pablo Ar\'es Gastesi}
\address{The Institute of Mathematical Sciences, Madras, India
and School of Maths, Tata Institute, Bombay, India}
\curraddr{School of Maths, Tata Institute, Bombay 400 005, India}
\email{pablo@@math.tifr.res.in}
\date{\today}
\vspace*{5cm}\maketitle
\begin{abstract}
The Torelli group of a compact non-orientable Klein surface is
the subgroup of the modular group consisting of the mapping classes that
act trivially on the first homology group of the surface. We
prove that if a surface has genus at least $3$, then the Torelli group
acts fixed points free on the \T space of the surface.
That gives an embedding of
the Torelli space of a Klein surface in the Torelli space of its complex
double. We also construct real tori associated to Klein surfaces, which we
call the Jacobian of the surface. We prove that this Jacobian is
isomorphic to a component of the real part of the Jacobian of the complex
double.
\end{abstract}

\section{Statement of results}

Klein surfaces are the natural generalization of Riemann surfaces to the
non-orientable situation: one considers holomorphic and anti-holomorphic
changes of coordinates. One of the points of interest in the
study of Klein surfaces is to determine which results of the theory of
deformation of Riemann surfaces hold for the non-orientable
case. A common approach to this problem is to consider a Klein
surface $\Sigma$, as a Riemann surface $\Sigma^c$, with an
anti-holomorphic involution
$\sigma$, and thus one wants to find $\sigma$-invariant objects. In this
paper we follow these two points of view to
show two related results in the theory of Klein surfaces,
that is, we will do some constructions on Klein surfaces, and then find
the corresponding invariant objects related to the Riemann surface
$\Sigma^c$. More precisely, we
construct the {\bf Torelli space} $Tor(\Sigma)$, and prove that it can
be identified with the set of fixed points of an involution on
$Tor(\Sigma^c)$. We also construct the {\bf Jacobian variety} $J(\Sigma)$
of $\Sigma$ by integrating a basis of the space of real
harmonic forms over the free part of
${\mathrm H}_1(\Sigma,\Z)$. We prove that $J(\Sigma)$ is
isomorphic to a component of the real part of the Jacobian $J(\Sigma^c)$
of the complex double $\Sigma^c$ of $\Sigma$.

Given a compact smooth non-orientable surface $\Sigma$, the {\T space}
$T(\Sigma)$ of $\Sigma$ is defined as $T(\Sigma)={\cal
M}(\Sigma)/Diff_0(\Sigma)$, where ${\cal M}(\Sigma)$ is the set of Klein
surface structures on $\Sigma$ that agree with the given smooth
structure, and $Diff_0(\Sigma)$ is the group of diffeomorphisms of
$\Sigma$ homotopic to the identity \cite[pg. 145]{sep:book}.
We will use $\Sigma$ for a Klein
surface, if it is clear from the context what the structure is, or we
will write $(\Sigma,X)$ if we need to specify more.
The {\bf modular} or {\bf mapping
class group}, $Mod(\Sigma)=Diff(\Sigma)/Diff_0(\Sigma)$, acts on
$T(\Sigma)$ by pull-back of dianalytic structures (see \S $2$).
The {\bf Torelli group} $U(\Sigma)$ is the subgroup of
$Mod(\Sigma)$ consisting of the mapping classes that act trivially on
${\mathrm H}_1(\Sigma,\Z)$. The parallel result to the following theorem
is a classical fact on Riemann surfaces.
\setcounter{section}{3}
\setcounter{thm}{0}
\begin{thm}Let $\Sigma$ be a compact non-orientable surface of genus $g
\geq 3$. Let $[f]\in Mod(\Sigma)$, and suppose that there exists a Klein
surface structure $X$ on $\Sigma$ such that $f:(\Sigma,X)\rightarrow
(\Sigma,X)$ is dianalytic. Then $[f]=[id]$. Therefore, $U(\Sigma)$ acts
fixed-points free on $T(\Sigma)$, and the Torelli space
$Tor(\Sigma)=T(\Sigma)/U(\Sigma)$ is a smooth real manifold of dimension
$3g-6$.
\end{thm}
Assume now that $\Sigma$ has a fixed Klein surface structure.
Then there exists an
unramified double covering of $\Sigma$ by a Riemann surface $\Sigma^c$,
known as the {\bf complex double}. Moreover, $\Sigma$ is isomorphic to
$\Sigma^c/<\!\sigma\!>$, where $\sigma$ is an anti-holomorphic
involution. The mapping $\sigma$ induces involutions
$\sigma^*$ and $\tilde\sigma$ on $T(\Sigma)$ and $Tor(\Sigma)$,
respectively. It is a well known fact that $T(\Sigma)$ can be identified
with the set of fixed points of $\sigma^*$. A similar result holds for
Torelli spaces, as the next proposition shows.
\setcounter{thm}{2}
\begin{prop}The Torelli space $Tor(\Sigma)$ can be identified with the
set of fixed points of $\tilde\sigma$ on $Tor(\Sigma^c)$.
\end{prop}

Torelli spaces are intimately related to the Jacobian
variety of a compact Riemann surface.
Recall that this variety $J(\Sigma^c)$, is a $g$-dimensional complex
torus ($g$ is the genus of $\Sigma^c$) given by $\C^g/\Gamma$, where
$\Gamma$ is the lattice generated by integration of a basis of holomorphic
forms on $\Sigma^c$ over a basis of ${\mathrm H}_1(\Sigma^c,\Z)$.
We can also construct the Jacobian by
considering the lattice $\Gamma'$ generated by integration of harmonic
forms and then taking the quotient
$\R^{2g}/\Gamma'$, which is a real torus. The Hodge-$*$ operator gives
a complex structure to this real torus in such a way that it becomes
$J(\Sigma^c)$. This point of view can be generalized to construct a
Jacobian variety $J(\Sigma)$ of a non-orientable Klein surface.
\setcounter{section}{4}
\setcounter{thm}{0}
\begin{thm}
Let $\Sigma$ be a compact non-orientable surface of
genus $g\geq3$. Then we can associate to $\Sigma$ a real torus of
dimension $g-1$, the {\bf Jacobian variety} $J(\Sigma)$ of $\Sigma$,
such that $J(\Sigma)$ is isomorphic to any component of the real part of
the Jacobian $J(\Sigma^c)$ of the complex double. This last set is
defined as the set of fixed points of the symmetry $\sigma_1$ of
$J(\Sigma^c)$ induced by $\sigma$.
\end{thm}
{\bf Acknowledgments}: the idea of using harmonic forms to construct
$J(\Sigma)$ was suggested by S. Nag; I would like to thank him for
useful conversations regarding this topic.
I would like also to express my gratitude to D. S.
Nagaraj and R. R. Simha for many helpful conversations while this paper
was being written.
\setcounter{section}{1}
\setcounter{thm}{0}
\section{Some general facts about Klein surfaces}

A {\bf Klein surface} (or {\bf dianalytic}) structure $X$ on a surface
without boundary $\Sigma$ is a covering by open sets $U_i$, and a
collection of homeomorphisms $z_i:U_i\rightarrow V_i$, where
$V_i\subset\C$ are open sets, such that $z_i\circ z_j^{-1}$ is
holomorphic or anti-holomorphic, whenever $U_i\cap U_j\neq\emptyset$
\cite{all:klein}. Observe that a
Klein surface structure on an orientable surface is just a pair of
conjugate Riemann surface structures \cite{natan:klein}.

A compact non-orientable surface $\Sigma$ is homeomorphic to
the connected sum of $g\geq 1$ real projective planes \cite{blackett:topo}.
The number $g$ is
called the {\bf genus} of $\Sigma$. If $g=2n+1$, then the fundamental
group of $\Sigma$ has a presentation given by generators $c$,
$a_1,\ldots,a_n$, $b_1,\ldots,b_n$, satisfying
$c^2\prod_{j=1}^n[a_j,b_j]=1$, where $[a,b]=aba^{-1}b^{-1}$. If the
genus is even, $g=2n+2$, then we can choose generators $c$, $d,$
$a_1,\ldots,a_n$, $b_1,\ldots,b_n$, satisfying the relation
$c^2d^2\prod_{j=1}^n[a_j,b_j]=1$. An alternative presentation for this
latter case is given by generators $\gamma$, $\delta$, $a_1,\ldots,a_n,$
$b_1,\ldots,b_n,$ and the relation
$\gamma\delta\gamma^{-1}\delta\prod_{j=1}^n[a_j,b_j]=1$.
For the rest of this paper,
we will assume that all surfaces are compact without boundary. We will
further assume that non-orientable surfaces have genus $g\geq 3$, while
orientable surfaces satisfy $g\geq 2$.

The {\bf complex double} \cite{all:klein} of a Klein surface $\Sigma$
of genus $g$ is a triple $(\Sigma^c,\pi,\sigma)$, where:\\
\noindent (1) $\Sigma^c$ is a Riemann surface of genus $g-1$;\\
\noindent (2) $\pi:\Sigma^c\rightarrow\Sigma$ is an unramified double
covering;\\
\noindent (3) there exist local coordinates $z$ and $w$ on $\Sigma^c$
and $\Sigma$, respectively, such that the function $w\circ\pi\circ
z^{-1}$ is either holomorphic or anti-holomorphic (i.e. $\pi$ is a
morphism of Klein surfaces);\\
\noindent (4) $\sigma:\Sigma^c\rightarrow\Sigma^c$ is a
symmetry such that $\pi\circ\sigma=\pi$.

Let $S$ be a compact orientable surface, with a fixed orientation and a
smooth structure. The {\bf \T space} $T(S)$ of $S$ is $T(S)={\cal
M}(S)/Diff_0(S)$, where ${\cal M}(S)$ is the set of Riemann surface
structures on $S$ that agree with the given orientation and smooth
structure \cite{sep:book}.
The classical definition of $T(S)$ involves quasiconformal
mappings; to see that it is equivalent to the above definition, it
suffices to observe that on a compact surface, any homeomorphism is
homotopic to a smooth one, and diffeomorphisms are quasiconformal.
The {\bf modular} or {\bf mapping class} group $Mod(S)$
is the group of homotopic classes of orientation preserving
diffeomorphisms of $S$, that is $Mod(S)=Diff^+(S)/Diff_0(S)$. This group
acts on $T(S)$ by pull-back of complex structures: if $[f]\in Mod(S)$,
and $[X]\in T(S)$, then $[f]^*([X])=[f^*(X)]$, where $f^*(X)$ is the
Riemann surface structure on $S$ that makes
$f:(\Sigma,f^*(X))\rightarrow (\Sigma,X)$ biholomorphic. However, this
action has fixed points; it is therefore interesting to find subgroups of
$Mod(\Sigma)$ that act without fixed points on $T(S)$. A subgroup $G$ of
$Mod(S)$ has the {\bf Hurwitz-Serre} property \cite{nag:teic}
if $G$ satisfies that for
any element $[g]\in G$ such that there exists an $[X]\in{\cal M}(S)$
with $g:(S,X)\rightarrow (S,X)$ biholomorphic, one has that $[g]=[id]$.
A group with this property will act fixed-points free on \T space and,
therefore, the quotient $T(S)/G$ will be a smooth finite dimensional
complex manifold. The {\bf Torelli group} $U(S)=\{[f]\in Mod(S);~
f~\mathrm{acts~trivially~on~H}_1(S,\Z)\}$, is known to satisfy the
Hurwitz-Serre property \cite{fk:book}.
The quotient space $Tor(S)=T(S)/U(S)$ is called
the {\bf Torelli space} of $S$.

The Jacobian variety $J(S)$ of a compact Riemann surface is an abelian
variety constructed as follows: let ${\cal
B}^c=\{\alpha_1,\ldots,\alpha_g,\beta_1,\ldots,\beta_g\}$ be a symplectic
basis of ${\mathrm H}_1(S,\Z)$. Then we can find a dual basis for
${\mathrm H}^0(S,\Omega_S^1)$, the space of holomorphic forms on $S$,
consisting of forms $\{\omega_1,\ldots,\omega_g\}$, satisfying
$$\int_{\alpha_j}\omega_k=\begin{cases}
 1 & \text{if } j=k,\\
 0 & \text{otherwise}.\end{cases}$$
Let $\Gamma^c$ be the lattice
on $\C^g$ generated by the vectors $(\int_c\omega_1,\ldots,\int_c\omega_g)$,
$c\in{\cal B}^c$; then we define $J(S)=\C^g/\Gamma^c$.
\section{Torelli groups of non-orientable compact surfaces}

In this section we will show that the Torelli group of a compact
non-orientable Klein surface $\Sigma$ has the Hurwitz-Serre property.
The main idea of the proof is to see that, if a diffeomorphism of a
Klein surface $\Sigma$ acts trivially on ${\mathrm H}_1(\Sigma,\Z)$,
its orientation preserving lift to the complex double will act
trivially on ${\mathrm H}_1(\Sigma^c,\Z)$; then we use the fact that the
Hurwitz-Serre property is satisfied for Riemann surfaces.
The quotient space $Tor(\Sigma)=T(\Sigma)/U(\Sigma)$ is a
smooth real manifold. We will show that $Tor(\Sigma)$ can be
identified with the set of fixed points of a symmetry $\tilde\sigma$ on
$Tor(\Sigma^c)$.

Let us start with a smooth non-orientable surface $\Sigma$, of genus
$g=2n+1$, and a diffeomorphism $f:\Sigma\rightarrow\Sigma$ such that $f$
acts trivially on ${\mathrm H}_1(\Sigma,\Z)$. We can find a unique
orientation preserving diffeomorphism $\tilde f$ of $\Sigma^c$ such that
the following diagram commutes \cite{sep:spaces}:
$$\diagram
\Sigma^c\rto^{\tilde f}\dto_\pi 	&  \Sigma^c\dto^\pi \\
\Sigma\rto^f				&  \Sigma.\enddiagram$$
We want to show that the mapping ${\tilde f}_\#$ induced by $\tilde f$
on ${\mathrm H}_1(\Sigma^c,\Z)$ is trivial. For that purpose we need
to recall the way $\Sigma^c$ is constructed, from the
topological viewpoint. The reader can find more details in
\cite{blackett:topo}.

By the presentation of the fundamental group of $\Sigma$, we can identify
this surface with a $(4n+2)$-polygon, whose sides are labeled
to satisfy the relation of the fundamental group. Then $\Sigma^c$
is given by two polygons with boundary relations:
$$c_1c_2\prod_{j=1}^n[a_{j,1},b_{j,1}]=1\hspace{5mm}
\mathrm{and}\hspace{5mm}
c_2c_1\prod_{j=1}^n[a_{j,2},b_{j,2}]=1.$$
To obtain a single relation, we find the value of $c_2$ on the right
hand side equation and substitute it on the left hand one
(equivalently, we glue the two polygons by the $c_2$ sides):
$$c_2=(\prod_{j=1}^n[b_{n+1-j,2},a_{n+1-j,2}])c_1^{-1};$$
therefore
$$c_1\big(\prod_{j=1}^n[b_{n+1-j,2},a_{n+1-j,2}]\big)c_1^{-1}
\big(\prod_{j=1}^n[a_{j,1},b_{j,1}]\big)=$$
$$\big(\prod_{j=1}^n[c_1b_{n+1-j,2}c_1^{-1},c_1a_{n+1-j,2}c_1^{-1}]\big)
\big(\prod_{j=1}^n[a_{j,1},b_{j,1}]\big)=1.$$
{}From this formula we see that $\Sigma^c$ is a compact surface of genus
$g-1=2n$; we can choose the following paths as generators of the
fundamental group of $\Sigma$:
$$\alpha_1=c_1b_{n,2}c_1^{-1},\ldots,\alpha_n=c_1b_{n,2}c_1^{-1},
\alpha_{n+1}=a_{1,1}\ldots,\alpha_{2n}=a_{n,1},$$
$$\beta=c_1a_{n,2}c_1^{-1},\ldots,\beta_n=c_1a_{n,2}c_1^{-1},
\beta_{n+1}=b_{1,1}\ldots,\beta_{2n}=b_{n,1}.$$
These loops satisfy $\prod_{j=1}^{2n}[a_j,b_j]=1$. Let $\cal B$ and
${\cal B}^c$ denote the basis of ${\mathrm H}_1(\Sigma,\Z)$ and
${\mathrm H}_1(\Sigma^c,\Z)$ induced by the two given sets of generators of the
corresponding fundamental groups. By an abuse of notation, we will use
the same letters for the elements of the fundamental group and their
classes in homology. We can see that ${\cal B}^c$ is a symplectic basis
of ${\mathrm H}_1(\Sigma^c,\Z)$; that is, the intersection matrix is
given by
$J=\left(\begin{matrix}
0 & {\mathrm I} \\
-{\mathrm I} & 0 \end{matrix}\right) ,$ where ${\mathrm I}$ is the identity
matrix. The covering map $\pi:\Sigma^c\rightarrow\Sigma$ induces a
mapping on homology, with associated matrix
$$\pi_\#=\left(\begin{matrix}
0  & 0  & 0  & 0  \\
0  & {\mathrm I}  & K  & 0 \\
K  & 0  & 0  & {\mathrm I} \end{matrix}\right),$$ with respect to ${\cal
B}^c$ and $\cal B$. The matrix $K$ is given by
$$K=\left(\begin{matrix}
0      & \cdots & \cdots & 1      \\
0      & \cdots &    1   & 0      \\
\vdots &        &        & \vdots \\
1      &   0    & \cdots &    0    \end{matrix}\right).$$
The symmetry $\sigma$ maps $a_{j1}$ (resp. $b_{j1}$) to $a_{j2}$ (resp.
$b_{j2}$); it is not difficult to see that the map $\sigma_\#$
induced on ${\mathrm H}_1(\Sigma^c,\Z)$ is given by $\sigma_\#=K$.
Let ${\tilde f}_\#=\left(A_{jk}\right)_{j,k=1}^4$;
since $f_\#$ acts trivially on ${\mathrm
H}_1(\Sigma,\Z)$, we have $\pi_\#{\tilde f}_\#=\pi_\#$. By the
uniqueness of ${\tilde f}_\#$ we get
${\tilde f}_\#\sigma_\#=\sigma_\#{\tilde f}_\#$.
Finally, ${\tilde f}^t_\#J{\tilde
f}_\#=J$, where ${\tilde f}_\#^t$ is the transpose of ${\tilde f}_\#$,
since $\tilde f$ preserves the intersection matrix
(\cite[theorem N13, pg. 178]{mag:comb}).
The condition $\pi_\#{\tilde f}_\# =
\pi_\#$ is equivalent to the following set of equations:
$$\left\{\begin{array}{lclccclcl}
A_{21} + KA_{31} & = & 0 & & & & KA_{11} + A_{41} & = & K \\
A_{22} + KA_{32} & = & {\mathrm I} & & & & KA_{12} + A_{42} & = & 0 \\
A_{23} + KA_{33} & = & K & & & & KA_{13} + A_{43} & = & 0 \\
A_{24} + KA_{34} & = & 0 & & & & KA_{14} + A_{44} & = & {\mathrm I}.
\end{array} \right .$$
Therefore, the matrix ${\tilde f}_\#$ can be written as
$${\tilde f}_\# = \left(\begin{matrix}
A_{11} & A_{12} & A_{13} & A_{14} \\
A_{21} & A_{22} & A_{23} & A_{24} \\
-KA_{21} & K-KA_{22} & {\mathrm I}-KA_{23} & -KA_{24} \\
K-KA_{11} & -KA_{12} & -KA_{13} & {\mathrm I}-KA_{14} \end{matrix}\right).$$
We now use the fact that $\tilde f$ and $\sigma$ commute, to obtain the
following relations among the entries of the matrix ${\tilde f}_\#$:
\begin{equation}\left\{ \begin{array}{lclccclcl}
A_{14}K & = & K(K-KA_{11}) & & & & A_{24}K & = & K(-KA_{21}) \\
A_{13}K & = & K(-KA_{12}) & & & &  A_{23}K & = & K(K-KA_{22}) \\
A_{12}K & = & K(KA_{13}) & & & &  A_{22}K & = & K({\mathrm I}-KA_{23}) \\
A_{11}K & = & K({\mathrm I}-KA_{14})& & & &  A_{21}K & = & K(-KA_{24}).
\end{array}\right .\label{eq:f}\end{equation}
Consider the equation ${\tilde f}_\#^t J {\tilde f}_\# = J$; looking
at the first row of the matrices on both sides of the equality, we get,
after using \eqref{eq:f} to simplify the result,
$$\left\{\begin{array}{rcr}
A_{21}^tK-KA_{21} & = & 0 \\
A_{11}^tK-KA_{22} & = & 0 \\
-2{\mathrm I}+A_{11}^t+KA_{22}K & = & 0 \\
A_{21}^t+KA_{21}K & = & 0
\end{array}\right .$$
Solving these equations, we get $A_{21}=0$ and $A_{11}={\mathrm I}$, which
imply that $A_{24}=0$ and $A_{14}=0$.
Using this, we now consider the equality between the second rows of the
matrices ${\tilde f}_\#^t J {\tilde f}_\#$ and $J$, to obtain
$A_{22}={\mathrm I}$, and $A_{12}=0$. By \eqref{eq:f}, we get $A_{23}=0$
and $A_{13}=0$. Therefore, we have that ${\tilde f}_\#$ is the identity
matrix.

{\bf Remark:} if we would have chosen the orientation reversing lift of $f$,
say ${\tilde f}_1$, then ${\tilde f}_1={\tilde f}\sigma$, so $({\tilde
f}_1)_\#={\tilde f}_\#\sigma_\#=\sigma_\#$.

If $\Sigma$ has even genus $g=2n+2$, we use the first of the two presentations
of its fundamental group given in \S $2$.
We have that $\Sigma^c$ is given by two polygons with boundary relations:
$$c_1c_2d_1d_2\prod_{j=1}^n[a_{j,1},b_{j,1}] = 1\hspace{5mm}
\mathrm{and}\hspace{5mm}
c_2c_1d_2d_1\prod_{j=1}^n[a_{j,2},b_{j,2}] = 1.$$
{}From the second equation we get
$$d_2=c_1^{-1}c_2^{-1}\big(\prod_{j=1}^n[b_{n+1-j,2},a_{n+1-j,2}]\big)
d_1^{-1},$$
which reduces the first equation to
$$c_1c_2d_1c_1^{-1}c_2^{-1}\big(\prod_{j=1}^n[b_{n+1-j,2},a_{n+1-j,2}]\big)
d_1^{-1}\big(\prod_{j=1}^n[a_{j,1},b_{j,1}]\big) = $$
$$c_1c_2d_1c_1^{-1}c_2^{-1}d_1^{-1}
\big(\prod_{j=1}^n[d_1b_{n+1-j,2}d_1^{-1},d_1a_{n+1-j,2}d_1^{-1}]\big)
\big(\prod_{j=1}^n[a_{j,1},b_{j,1}]\big) = 1.$$
We therefore obtain that the fundamental group of $\Sigma^c$ is
generated by the loops
$$\alpha_1=c_1d_1^{-1},~\alpha_2=d_1b_{n,2}d_1^{-1},\ldots,
\alpha_{n+1}=d_1b_{1,2}d_1^{-1},
\alpha_{n+2}=a_{1,1},\ldots,\alpha_{2n+1}=a_{n,1},$$
$$\beta_1=d_1c_2,~\beta_2=d_1a_{n,2}d_1^{-1},\ldots,
\beta_{n+1}=d_1\a_{1,2}d_1^{-1},
\beta{n+2}=a_{1,1},\ldots,\beta_{2n+1}=a_{n,1},$$
satisfying the relation
$\prod_{j=1}^{2n+1}[\alpha_j,\beta_j]=1.$
The basis $\{\alpha_1, \alpha_2, \ldots, \alpha_{2n+1}, \beta_1, \beta_2,
\ldots,$\newline
\noindent $\beta_{2n+1} \}$
is symplectic, but computations are easier if we
rearrange the basis as
${\cal B}^c=\{ \alpha_1, \beta_1, \alpha_2, \ldots,\alpha_{2n+1},
\beta_2, \ldots, \beta_{2n+1} \}$,
whose intersection matrix is given by:
$$J=\left(\begin{matrix}
 N &  0 &  0 &  0 &  0 & \\
 0 &  0 &  0 &  {\mathrm I} &  0 & \\
 0 &  0 &  0 &  0 &  {\mathrm I} & \\
-{\mathrm I} &  0 &  0 &  0 &  0 & \\
 0 & -{\mathrm I} &  0 &  0 &  0 & \end{matrix}\right) ,$$
where $N=\left(\begin{matrix} 0 & 1 \\ -1 & 0 \end{matrix}\right) .$
With respect to $\cal B^c$, the symmetry $\sigma$ has the following
matrix representation, for its action on ${\mathrm H}_1(\Sigma^c,\Z)$:
$$\sigma_\#=\left(\begin{matrix}
M & 0 & 0 & 0 & 0 \\
0 & 0 & 0 & K & 0 \\
0 & 0 & K & 0 & 0 \\
0 & K & 0 & 0 & 0 \\
K & 0 & 0 & 0 & 0 \end{matrix}\right) ,$$
where the matrix
$M=\left(\begin{matrix} 1& 0 \\
2 & -1 \end{matrix}\right) .$
To obtain the matrix
$M$, observe that in homology one has $\alpha=c_1+d_1$ and
$\beta=d_1+c_2$, so $\sigma(\alpha)=c_2+d_2$,
$\sigma(\beta)=d_2+c_1.$ Substituting the value of $d_2$ obtained
previously, we get $M$.
Similarly we have that the expression for the action of the covering
map $\pi$ on homology is given by:
$$\pi_\#=\left(\begin{matrix}
{\mathrm I} & 0 & 0 & 0 & 0 \\
0 & 0 & {\mathrm I} & K & 0 \\
0 & K & 0 & 0 & {\mathrm I} \end{matrix}\right) .$$
Computing in a  way similar to the odd genus case, we get that if
$f:\Sigma\rightarrow\Sigma$ is a diffeomorphism of $\Sigma$ that acts
trivially on ${\mathrm H}_1(\Sigma,\Z)$, and $\tilde f$ is a lift of $f$
to $\Sigma^c$, then the action of this last mapping on ${\mathrm
H}_1(\Sigma^c,\Z)$ is given by either the identity or $\sigma_\#$,
depending on whether $\tilde f$ is orientation preserving or reversing.

Recall that a map $f:\Sigma\rightarrow\Sigma$ on a Klein
surface is {\bf dianalytic} if, when expressed in local coordinates $(U,z)$,
$f\circ z^{-1}$ is either holomorphic or anti-holomorphic. The above
computations show that $U(\Sigma)$ satisfies the equivalent of the
Hurwitz-Serre property.
\begin{thm}Let $\Sigma$ be a compact non-orientable surface of genus
$g\geq 2$. Let $[f]\in Mod(\Sigma)$, and suppose that there exists a
Klein surface structure $X$ on $\Sigma$ such that
$f:(\Sigma,X)\rightarrow (\Sigma,X)$ is dianalytic. Then $f$ is
homotopic to the identity.
\end{thm}
\begin{cor}The Torelli space $Tor(\Sigma)=T(\Sigma)/U(\Sigma)$
is a smooth manifold of real dimension $3g-6$.
\end{cor}
\begin{proof}{of the Proposition} Since $f$ is dianalytic on the Klein surface
$(\Sigma,X)$, the orientation preserving lift $\tilde f$ is
biholomorphic on the Riemann surface $\Sigma^c$. But
then, since the genus of $\Sigma^c$ is at least $2$, we have that
$\tilde f$ is the identity, which proves the result.
\end{proof}
The involution $\sigma$ of $\Sigma^c$ induces a symmetry $\sigma^*$ on
$T(\sigma^c)$. The \T space $T(\Sigma)$ can be identified with the set
of fixed points of $\sigma^*$, which proves the corollary.
It is clear that $\sigma^*$ descends to a symmetry
$\tilde\sigma$ on $Tor(\Sigma^c)$.
\begin{prop}The Torelli space $Tor(\Sigma)=T(\Sigma)/U(\Sigma)$ can be
identified with the set of fixed points of $\tilde\sigma$ in
$Tor(\Sigma^c)$.\end{prop}
\begin{pf}The proof follows immediately from the definition of Torelli
spaces. In fact, we have that two elements $[X_1]$ and $[X_2]$ of ${\cal
M}(\Sigma)$ project to the same point in $Tor(\Sigma)$ {\it if and only if}
there exists a diffeomorphism $h\in Diff(\Sigma)$ such that
$h_\#:{\mathrm H}_1(\Sigma,\Z)\rightarrow {\mathrm H}_1(\Sigma,\Z)$ is
the identity, and $h:(\Sigma,X_2)\rightarrow (\Sigma,X_1)$ is
dianalytic. The rest of the proof is similar to the proof that
$T(\Sigma)$ can be identified with the set of fixed points of
$\sigma^*$; see \cite{sep:book} for more details.\end{pf}
\section{Jacobi varieties of Klein surfaces}

Throughout this section, $\Sigma$ will denote a fixed compact
non-orientable Klein surface of genus $g\geq 3$, and $\Sigma^c$ its
complex double. We can take $\Sigma^c$ to be defined by a
polynomial, $p(z,w)=0$, with real coefficients (\cite{all:klein} and
\cite{natan:gordon}). Then the involution
$\sigma$ is given by $\sigma(z,w)=(\overline z,\overline w)$, and
conjugation $z\mapsto\overline z$ in $\C^{g-1}$ induces an involution
$\sigma_1$ on the Jacobian $J(\Sigma^c)$. The set of fixed points of
$\sigma_1$, that is, the {\bf real part} of $J(\Sigma^c)$, is a real manifold
of dimension $g-1$; the pair $J(\Sigma^c,\sigma_1)$ is usually
considered as the Jacobian of $\Sigma$.
On a Klein surface the concept of harmonic forms makes
sense; it is not difficult to see that
the space $\HS$ of such forms has dimension precisely $g-1$. One
can choose a basis of ${\mathrm H}_1(\Sigma,\Z)_f$, the free part of
${\mathrm H}_1(\Sigma,\Z)\cong\Z^{g-1}\oplus\Z/2Z$,
and a dual basis for $\HS$;
these two basis generate a lattice $\Gamma$ in
$\R^{g-1}$. We will call the real torus $\R^{g-1}/\Gamma$
the {\bf Jacobian variety} of $\Sigma$, and denote it by $J(\Sigma)$. On the
other hand, the real part of a holomorphic form (on $\Sigma^c$) is a
harmonic form, so one can expect some relationship between $J(\Sigma)$
and the real part of $J(\Sigma^c)$. We prove that, in fact, $J(\Sigma)$
is isomorphic to a component of the set of fixed points of $\sigma_1$ in
$J(\Sigma^c)$.

A continuous function $f:W\rightarrow\R$ defined on an open set of a
Klein or Riemann surface is called {\bf harmonic} if for any local
coordinate $(U,z)$, with $U\cap W\neq\emptyset$, the function $f\circ
z^{-1}$ is harmonic. Since precomposition with holomorphic and
anti-holomorphic functions preserves harmonicity, the above definition
makes sense. Actually, a Klein surface is the most general surface in
which the notion of harmonic function is well defined \cite{all:klein}.
Similarly, a (real) form $\psi$ is {\bf harmonic} if it
can be written locally as $\psi=df$, where $f$ is harmonic. We will
denote by $\HS$ the space of harmonic forms on $\Sigma$. Let $\sigma^*$
be the pull-back map induced by $\sigma$ on forms on $\Sigma^c$.
If $\omega=gdz$, with $g$ holomorphic, we have that $\sigma^*(\omega) =
g(\sigma)\sigma_{\overline z}d\overline z$, so $\sigma^*$ is
anti-holomorphic. We also have that, for any holomorphic form $\omega\in
{\mathrm H}^0(\Sigma^c,\Omega^1)$, and for any cycle $c$ on $\Sigma^c$,
$\int_c\sigma^*(\omega)=\int_{\sigma(c)}\omega$. Observe that this last
equality agrees with \cite{natan:gordon}, while it differs
of \cite{silhol:comess} and \cite{all:klein}, since these two authors
define $\sigma^*$ as the conjugate of our definition (in order
to have that $\sigma^*$ preserves holomorphic forms).

To compute the dimension of $\HS$, it suffices to observe that
$\sigma^*$ takes harmonic forms to harmonic forms; therefore, $\HS$ will
be isomorphic to the set of fixed points of $\sigma^*$ in $\HSc$. By
Hodge theory, $\sigma^*$ acts like $\sigma_\#=K$; so $dim\,\HS=g-1$.
This result agrees with \cite{all:obst}.

In order to justify later computations, we need the following lemma.
\begin{lemma}[Duality Lemma]On a non-orientable compact Klein surface
$\Sigma$, of genus $g \geq 3$, the space of harmonic forms, $\HS$,
and the dual space to the homology with real coefficients,
${\mathrm H}_1(\Sigma,\R)^*$, are isomorphic.\end{lemma}
\begin{pf}From Differential Topology \cite[Theorem 15.8]{bott:forms}
we know that the \v{C}ech cohomology with coefficients in the constant
presheaf $\Bbb Z$, ${\mathrm H}^1_\Z(\Sigma)^*,$ is isomorphic to the
singular cohomology ${\mathrm H}^1(\Sigma)$. Furthermore, by the Universal
Coefficients Theorem
\cite[Corollary 15.14.1]{bott:forms},
we have that the space ${\mathrm H}^1(\Sigma)$ is isomorphic to the free
part of ${\mathrm H}_1(\Sigma)$, which is just ${\Bbb Z}^{g-1}$.
Tensoring with $\R$ we
have that ${\mathrm H}_\R^1(\Sigma)$ is isomorphic to $\R^{g-1}$.
On the other hand, ${\mathrm H}_\R^1(\Sigma)$ is isomorphic to the de Rham
cohomology ${\mathrm H}_{DR}^1(\Sigma)$ \cite[8.9, 9.8 and 14.28]{bott:forms}.
By the compactness of $\Sigma$, on each de Rham class there exists at most a
harmonic form. A counting of dimensions shows that there exists
exactly one harmonic form, which completes the proof. Therefore, $\HS$
and ${\mathrm H}_{DR}(\Sigma)$ are isomorphic, and since this last space is
isomorphic (by integration) to ${\mathrm H}_1(\Sigma,\R)$, we are done.
\end{pf}

To make matters more clear, we will construct $J(\Sigma)$ on the cases
of genus $3$ and $4$. The general case will follow easily from our
examples.  Let us start with a Klein surface $\Sigma$  of
genus $g=3$.
Then from \S $2$
we have that ${\cal B}=\{a,b\}$ is a basis of ${\mathrm
H}_1(\Sigma,\Z)_f$. The loops $\alpha_1$, $\alpha_2$, $\beta_1$,
$\beta_2$ of \S $3$ give a basis of ${\mathrm H}_1(\Sigma^c,\Z)$; but
for computational purposes, it is better to choose the basis
$${\cal B}^c=\{\gamma_1 = -(\alpha_2+\beta_1), \gamma_2 =
-(\alpha_1+\beta_2), \delta_1 = \alpha_1+\alpha_2+\beta_1, \delta_2 =
\alpha_1+\alpha_2+\beta_2\}.$$
Observe that the change of basis in ${\mathrm H}_1(\Sigma^c,\Z)$ is
given by the matrix
$$C=\left(\begin{matrix}
0 & -1 & 1 & 1 \\
-1 & 0 & 1 & 1 \\
-1 & 0 & 1 & 0 \\
0 & -1 & 0 & 1 \end{matrix}\right) ,$$
which has determinant equal to $1$, and satisfies $C^tJC=J$.
Therefore ${\cal B}^c$ is a symplectic basis. The mapping $\pi_\#$ gives
$\pi_\#(\gamma_1) = -2b,$ $\pi_\#(\gamma_2) = -2a.$ This suggests that
we should take a basis ${\cal B}_h=\{\phi_1, \phi_2\}$ of $\HS$ normalized
by $\int_a\phi_1=\int_b\phi_2=0$, and $\int_b\phi_1=\int_a\phi_2=-1/2.$
Observe that this normalization is possible because of the Duality
Lemma, and the fact that $a$ and $b$ are not torsion classes in
${\mathrm H}_1(\Sigma,\Z)$.
We can use the pull-back mapping $\pi^*$ induced by $\pi$ to get forms
$\psi_j=\pi^*(\phi_j)$ on $\Sigma^c$. These forms are real harmonic, so
$\omega_j=\psi_j+i*\phi_j$ are holomorphic forms (where $*$ stands for
the Hodge-$*$ operator). By the expression of
$\sigma$ and the formula relating $\sigma^*$ with integrals, we have
that
$$\int_{\gamma_j}\overline{\omega}_k = \int_{\gamma_j}\sigma^*(\omega_k)
= \int_{\sigma(\gamma_j)}\omega_k = \int_{\gamma_j}\omega_k .$$
In particular, we see that $\int_{\gamma_j}\omega_k$ is real. But since
by \cite[Theorem 1.0.7, pg. 74]{all:klein}
${\mathrm Re}\int_{\gamma_j}\omega_k = \int_{\gamma_j}\psi_k =
\int_{\pi(\gamma_j)}\phi_k,$ we have that ${\cal B}^*=\{\omega_1,
\omega_2\}$ is normalized with respect to ${\cal B}^c$. Let $P$ denote
the corresponding period matrix, that is the entries of this matrix are
given by $p_{jk} = \int_{\delta_k}\omega_j$. The mapping
$\sigma_\#:{\mathrm H}_1(\Sigma^c,\Z)\rightarrow {\mathrm
H}_1(\Sigma^c,\Z)$ has the following expression with respect to the
basis ${\cal B}^c$:
$$\sigma_\#=\left(\begin{matrix}
 1 & 0 & -2 & -1 \\
 0 & 1 & -1 & -2 \\
 0 & 0 & -1 &  0 \\
 0 & 0 &  0 & -1 \end{matrix}\right) .$$
This results agrees with the one obtained by Natanzon in
\cite{natan:gordon} except in that he gets all sign positive.
Nevertheless, the computation of the real part of $J(\Sigma^c)$ yields
the same result.
We have that $P$ satisfies $\overline{P}=A-P$, where
$A=\left(\begin{matrix}
 -2 & -1 \\
 -1 & -2 \end{matrix}\right) .$ The Jacobian variety of $\Sigma^c$ is
then given by $J(\Sigma^c) = \C^2/\Gamma^c$, where $\Gamma^c=\Z^2+P\Z^2$.
To compute the real part of $J(\Sigma^c)$, we write, for any $z\in\C^2$,
$z=P\alpha+\beta$, where $\alpha,~\beta\in\R^2$. Then
$\sigma_1(z)=\overline z\equiv z$ is equivalent to $P\alpha+\beta =
\overline{P}\alpha+\beta+Pn+m,$ for some $n,~m\in\Z^2$. The imaginary
part of this equation gives $({\mathrm Im} P)\alpha = -({\mathrm Im}
P)\alpha + ({\mathrm Im} P)n$. By Riemann bilinear relations
(\cite{fk:book}; see also \cite{simha:rs} for a nice introduction to
Jacobians and of Riemann surfaces from an algebro-geometric approach)
we have
that $({\mathrm Im} P)$ is invertible, so we obtain
$\alpha=n/2$. On the other hand, taking real parts in the above
equation we get $({\mathrm Re} P)\alpha+\beta = (A-{\mathrm Re}
P)\alpha+\beta+({\mathrm Re} P)n+m$. Since $\overline P = A-P$, we have
that $({\mathrm Re} P)=A-({\mathrm Re} P)$, so this equation reduces to
$0=\frac{1}{2}An+m$, or equivalently
$$\left\{\begin{array}{ccc}
n_1+\frac{n_2}{2}+m_1 & = & 0 \\
\frac{n_1}{2}+n_2+m_2 & = & 0, \end{array}\right .$$ where notation the
should be clear. This implies that $n_j\in 2\Z$, so $\alpha\in\Z^2$.
Therefore, the set of fixed points of $\sigma_1$ is given by the real
torus ${\mathrm Re}(J(\Sigma^c))=\{P\Z^2+\beta;~\beta\in\R^2\}/\Gamma^c$.
which agrees with the results obtained by Silhol
\cite[pgs. $349$ and 359]{silhol:comess} and Natanzon \cite{natan:klein}.
By the form of the lattice $\Gamma^c$,
it is clear that ${\mathrm Re}(J(\Sigma^c))\cong (\R/\Z)^2$.

In a similar way to the construction of $J(\Sigma^c)$, we can form a
lattice in $\R^2$ using the basis ${\cal B}_h=\{\phi_1,\phi_2\}$ and
${\cal B}$ of $\HS$ and
${\mathrm H}_1(\Sigma,\Z)_f$, respectively.
Let us denote this lattice by $\Gamma$. We define the {\bf
Jacobian variety} of $\Sigma$ as the quotient $J(\Sigma)=\R^2/\Gamma$.
It is clear that $[z]\mapsto[-\frac{1}{2}z]$ induces an isomorphism
between ${\mathrm Re}(J(\Sigma^c))$ and $J(\Sigma)$, which proves our
result for $g=3$. The general case of a surface of odd genus is done in
a similar way.

To see the even genus case, we take a surface with $g=4$,
and we choose the second
of the two presentations of the fundamental group of $\Sigma$ given in
\S $2$; i.e. the generators are the loops
$c$, $d$, $a$ and $b$, and the relation is
$cdc^{-1}d[ab] = 1$. To construct the
complex double we proceed as in \S $3$; we do not include the
computations here, since they are done as in \S $3$. We get that
the fundamental group of $\Sigma^c$ is generated by the loops
$$\alpha_1=c_2c_1, ~\alpha_2=(d_1^{-1}c_2)b_1(d_1^{-1}b_1c_2)^{-1},
{}~\alpha_3=a_2,$$
$$\beta_1=d_1^{-1}, ~\beta_2=(d_1^{-1}c_2)a_1(d_1^{-1}b_1c_2)^{-1},
{}~\beta_3=b_2,$$
We again change our basis of ${\mathrm H}_1(\Sigma^c,\Z)$ to
$${\cal B}^c=\{\gamma_1=\alpha_1,~ \gamma_2=-(\alpha_3+\beta_2),~
\gamma_3=-(\alpha_2+\beta_3),~ \delta_1=\alpha_1+\beta_1,~
\delta_2=\alpha_2+\alpha_3+\beta_2,~$$
$$\delta_3=\alpha_2+\alpha_3+\beta_3\}.$$
It is not hard to see that ${\cal B}^c$ is symplectic; one simply has
to chech that the matrix
$$C=\left(\begin{matrix}
 1 & 0  & 0  & 1 & 0 & 0 \\
 0 & 0  & -1 & 0 & 1 & 1 \\
 0 & -1 & 0  & 0 & 1 & 1 \\
 0 & 0  & 0  & 1 & 0 & 0 \\
 0 & 0  & -1 & 0 & 1 & 0 \\
 0 & -1 & 0  & 0 & 0 & 1 \end{matrix}\right) ,$$
satisfies $C^tJC=J$. The projection $\pi$ acts on these loops by
$\pi_\#(\gamma_1)=2c$, $\pi_\#(\gamma_2)=-2a$, $\pi_\#(\gamma_1)=-2b\,$;
so we take a basis $\{\phi_1, \phi_2, \phi_3\}$ of $\HS$, and normalize
it by requiring
$$\int_c\phi_1 = -\int_a\phi_2 = -\int_b\phi_3 = \frac{1}{2}, $$
and the other integrals equal to $0$.
As in the previous
situation, we have that if $\psi_j=\pi^*(\phi_j)$, then ${\cal B}^* = \{
\omega_j=\psi_j+i*\psi_j;~ j=1,2,3\}$ is a basis of holomorphic
$1$-forms dual to ${\cal B}^c$. It is not hard to see that $\sigma_\#$
is given by the following matrix, when computed with respect to ${\cal
B}^c$:
$$\sigma_\#=\left(\begin{matrix}
 1 & 0 & 0 & 2 & 0 & 0 \\
 0 & 1 & 0 & 0 & -2 & -1 \\
 0 & 0 & 1 & 0 & -1 & -2 \\
 0 & 0 & 0 & -1 & 0 & 0 \\
 0 & 0 & 0 & 0 & -1 & 0 \\
 0 & 0 & 0 & 0 & 0 & -1 \end{matrix}\right) ,$$
and the period matrix $P$ satisfies $\overline P=A-P$, where
$$A=\left(\begin{matrix}
 2 & 0 & 0 \\
 0 & -2 & -1 \\
 0 & -1 & -2 \end{matrix}\right).$$
The above matrix of $\sigma_\#$ is different from the one given in
\cite{natan:gordon}; we have not been able to obtain the matrix of that
reference, but nevertheless, we obtain similar result in the computation
of the real part of $J(\Sigma^c)$. In this case we have that the real part
of $J(\Sigma^c)$ (which is found in the same way that the $g=3$ case)
has two components, namely
$T_1=\{\Z^3n+\beta;~\beta\in\R^3,~ n=(n_1,n_2,n_3),~
n_1,n_2,n_3\in\Z\}$ and
$T_2=T_1+\Z^3(\frac{1}{2},0,0)^t$. We again obtained the results of
\cite{silhol:comess} and \cite{natan:klein}.
An isomorphism similar to the
previous case holds in this situation, except that we have $J(\Sigma)$
is isomorphic to any of the two sets $T_1$ or $T_2$.

{\bf Remark}: the results about the fixed points of $\sigma_1$ on
$J(\Sigma^c)$ can also be obtained from the expression of the period
matrix given in \cite[Proposition $4$, pg. 351]{silhol:comess}.

The above results can be put together in the following theorem:
\begin{thm}
Let $\Sigma$ be a compact non-orientable surface of
genus $g\geq3$. Then we can associate to $\Sigma$ a real torus of
dimension $g-1$, the {\bf Jacobian variety} $J(\Sigma)$ of $\Sigma$,
such that $J(\Sigma)$ is isomorphic to any component of the real part of
the Jacobian $J(\Sigma^c)$ of the complex double. This last set is
defined as the set of fixed points of the symmetry $\sigma_1$ of
$J(\Sigma^c)$ induced by $\sigma$.
\end{thm}

\end{document}